\documentclass[a4paper,reqno]{amsart}
\usepackage{pgf,tikz,pgfplots}
\usepackage{amsmath}
\usepackage{paralist}
\usepackage{graphics}
\usepackage{epsfig}
\usepackage{amsfonts}
\usepackage{amssymb}
\usepackage{hyperref,cite}

\usepackage{amsthm}
\usepackage{enumerate}
\usepackage[mathscr]{eucal}
\usepackage{eqlist}
\UseRawInputEncoding

\def\teich{{Teichm\"{u}ller}}
\def\mobi{{M{\"{o}}bius}}

\newtheorem{theorem}{Theorem}[section]
\newtheorem{lemma}[theorem]{Lemma}

\def\ifl{\iffalse }

\numberwithin{equation}{section}

\numberwithin{equation}{section}

\theoremstyle{remark}

\begin{document}
\baselineskip=1.2\baselineskip
\title[Strongly quasisymmetirc homeomorphisms]
{Strongly quasisymmetirc homeomorphisms  being compatible with  Fuchsian groups }

\author {Shengjin Huo}
\address{School of Mathematical Sciences, Tiangong University, Tianjin 300387, China} \email{huoshengjin@tiangong.edu.cn}
\author{Mengzhen Zhao}
\address{School of Mathematical Sciences, Tiangong University, Tianjin 300387, China} \email{2200152584@qq.com}
%
%


\thanks{This work was supported by the National Natural Science Foundation of China (Grant No.11401432).}
%
\begin{abstract}
 In this paper we first introduced a domain called generalized Dirichlet fundamental domain $\mathcal{F}^{*}$ for a Fuchsian group $G$ whose generators contain parabolic elements. This allows us to show that a quasisymmetric homeomorphism $h$ being compatible with a convergence Fuchsian group $G$ of first kind is a strongly quasisymmetric homeomorphism if and only if it has a quasiconformal extension  $f$ to the upper half plane $\mathbb{H}$ onto itself such that the induced measure $\lambda_{\mu}=|\mu|^{2}/Im(z)dxdy$ by the Beltrami coefficient $\mu$ of $f$ is a Carleson measure on the generalized Dirichlet fundamental domain $\mathcal{F}^{*}.$
 We also show that the above property  also holds for Carleson-Denjoy domains.
\end{abstract}

\maketitle
\section{ Introduction}
 The set $QS(\mathbb{R})$ of quasisymmetric homeomorphisms of the real axis $\mathbb{R}$ plays a very important role in the theory of Riemann surfaces, harmonic analysis, one dimensional real dynamical systems, and more generally, in the theory of {\teich} spaces. There are lots of papers investigated the properties of quasisymmetric homeomorphisms, see \cite{AZ}, \cite{Da}, \cite{GS},\cite{HS}, \cite{FHS},\cite{FKP}and the references of these papers therein.

 A subgroup $G$ of the group $PSL(2, \mathbb{R})$ (the group of all {\mobi} transformations of the upper half plane to itself) is called a Fuchsian group if it is a discrete subgroup of $PSL(2, \mathbb{R})$ with respect to the topology on $PSL(2, \mathbb{R})$ induced by the Euclidean distance on $\mathbb{R}^{4}$. For a Fuchsian group $G$, a quasisymmetric homeomorphism $h$ of the real axis to itself is called $G$-compatible if the conjugation $h\circ g\circ h^{-1}$ is a {\mobi} transformation for every $g\in G.$

We now present the context of the paper, see section 2 for the relevant definitions.
   An important problem of a long time has been to determine when a quasisymmetric homeomorphism is singular in the
usual sense: it maps a set with zero Lebesgue measure to the complement of zero
measure set (equal to full measure).  In 1979, Bowen \cite{Bo} proved that if $G$ is a cocompact Fuchsian group (the quotient Riemann surface $\mathbb{H}/G$ is compact, then all $G$-compatible quasisymmetric homeomorphisms except  {\mobi} transformations are always singular. Soon Sullivan \cite{Su1,Su2} extended Bowen's results to all cofinite groups.
   In \cite{BS}, Bishop and Steger got the following result: for a lattice group $G$ (i.e. G is finitely generated of first kind), there is a set $E\subset \mathbb{R}$ such that the Hausdorff dimensions of both $E$ and $h(\mathbb{R}\setminus E)$ are less than 1, where $h$ is a quasi-symmetric conjugating homeomorphism of the real axis $\mathbb{R}.$ Bishop and Steger's result implies that any conjugation $h$ of the real axis $\mathbb{R}$ of a lattice group $G$ must either be {\mobi} or strongly singular, i.e., $h$ maps a set of Hausdorff dimension $<1$ to the complement of a set of Hausdorff dimension $<1.$  In \cite{H}, we show that Bishop and Steger's result holds if and only if $G$ is geometric finite. If any quasi-symmetric homeomorphism $h$ conjugates a divergence Fuchsian group to another one, then $h$ is either {\mobi} or singular (see e.g., Agard \cite{A} , Tukia \cite{Tu} or Bishop \cite{Bi1}), i.e. $h$ is continuous but the derivation of $h$ vanishes almost everywhere in the real axis $\mathbb{R}$. In 1990, Astala and Zinsmeister \cite{AZ1} showed that the above singular property fails for the quasisymmetric conjugating homeomorphisms  of all convergence Fuchsian groups of the first kind.

On the contrary, it is also interesting to  determine when a  quasisymmetric homeomorphism being compatible with Fuchsian groups is absolutely continuous.
Carleson \cite{C} initiated such an investigation, giving a sufficient condition on the dilatation of a quasiconformal self-mapping of the unit disk to have an absolutely continuous boundary value. Along this direction, there are lots of researches about this problem, see \cite{AZ,BJ,CZ,Dy,SW}.

Since the singular properties of  quasisymmetric homeomorphisms which are  compatible with Fuchsian groups of divergence type are all singular.
The only interesting thing  is to consider  the  quasisymmetric homeomorphisms which are  compatible with Fuchsian groups of convergence type.
In this paper our motivation  is to continue investigating the absolutely continuous properties of  $h\in QS(\mathbb{R})$  which conjugates the
action of convergence type Fuchsian groups acting on the upper half-plane $\mathbb{H}$.

\section{ Terminologies and main results}
\subsection{ Quasisymmetric homeomorphisms.}
An orientation-preserving homeomorphism $f$
between planar domains $\Omega$ and $\Omega^{'}$ is called a quasiconformal mapping if
it satisfies the Beltrami equation
$$\bar{\partial} f=\mu\partial f,~~~~~a.e.~~ z\in \Omega,$$
where the measurable function  $\mu$ with essential norm $||\mu||_{\infty}<1$ is called the Beltrami coefficient of $f$.

 Let $\mathbb{H}$ be the upper half plane in the complex plane
$\mathbb{C}$,  $\mathbb{R}\cup\{\infty\}$ is its boundary at infinity.
 Let $M(\mathbb{H})$ denote the open unit ball of the Banach space $L^{\infty}(\mathbb{H})$
of all essentially bounded measurable functions on the upper half plane $\mathbb{H}.$    For every $\mu\in M(\mathbb{H})$, there is a quasiconformal self-homeomorphism $f$ of  the upper half plane $\mathbb{H}$ to itself with Beltrami coefficient $\mu$. Note that by the measurable Riemann mapping theorem we know that the quasiconformal self-homeomorphism $f$ is uniquely up to post-composition of elements of PSL(2, $\mathbb{R}$). We set the group of all orientation-preserving  quasiconformal self-homeomorphisms of the upper half plane $\mathbb{H}$ by $QC(\mathbb{H})$ and set the group of all increasing homeomorphism of $\mathbb{R}$ onto itself by $Hom^{+}(\mathbb{R})$.  An increasing homeomorphism $h\in Hom^{+}(\mathbb{R})$  is called quasisymmetirc if there is a constant $M \geq 1$ such that the $M$-condition,
 $$\displaystyle \frac{1}{M}\leq\displaystyle \frac{h(x+t)-h(x)}{h(x)-h(x-t)}\leq M,$$
 is satisfied for every symmetric triple $x-t$, $x$ and  $x+t$  in the real axis $\mathbb{R}$. Let $QS(\mathbb{R})$ be the group of all quasisymmetric self homeomorphisns of $\mathbb{R}.$

 Every $f\in QC(\mathbb{H}) $ can extend continuously to an increasing homeomorphism of $\mathbb{R}$ onto itself denoted by $h$. It is well known that $h$ is in $QS(\mathbb{R})$. A homeomorphism $h$ in $Hom^{+}(\mathbb{R})$ is said to be strongly quasisymmetric if there exist two positive constants $C_{1}$ and $C_{2}$ such that

 $$\displaystyle \frac{|h(E)|}{|h(I)|}\leq C_{1}(\displaystyle \frac{|E|}{|I|})^{C_{2}},$$
 where $I\subset \mathbb{R}$ is an interval, $E\subset I$ a measurable subset and $|\cdot|$ denotes the Lebesgue measure. Let $SQS(\mathbb{R})$ denote the set of all strongly quasisymmetric homeomorphisms of $\mathbb{R}$ onto itself.  A quasisymmetric homeomorphism $h$ is strongly quasisymmetric if and only if $h$ is locally absolutely continuous so that $h'$ belongs to the class of $A_{\infty}$ introduced by Muckenhoupt, see \cite{CF} or \cite{Ga}.

\subsection{ Fuchsian groups.} In this paper we call a {\mobi} group  $G$ Fuchsian group  if it acts on the upper half plane $\mathbb{H}$  properly discontinuously and freely, i.e a discrete subgroup of $ PSL(2, \mathbb{R})$.
Then the accumulation set of any  orbit $\{g(z); g\in G\}$ of $z\in \mathbb{H}$ does not depend on
$z$ and is called the limit set of $G$, denoted by $L(G)$.  The limit set $L(G)$ of $G$ is a closed subset of $\mathbb{R}.$
 If the limit set $L(G)$ agrees with $\mathbb{R}$, then $G$ is
said to be of the first kind. Otherwise, we say that $G$ is of the second kind.

A Fuchsian  group $G$ is of divergence type if $$\sum_{g\in G}\exp(-\rho(i,g(i)))=\infty,$$
where $\rho(i,g(i))$ is the hyperbolic distance between $i$ and $g(i)$. Otherwise, we say that it is of convergence type. All second kind groups are of convergence type. The points in $L(G)$ can be divided into different types by the ways in which they are approached by sequences in $Gz.$ A horizontal line or a circle tangent at the real line (without its point of tangency) is called a horocycle, and the tangency point of a horocycle is called its center. A disk with a horocycle boundary is called a horodisk.  A point $x$ in $L(G)$ is called horocyclic if for any $z$ in $\mathbb{H}$ its orbit $Gz$ meets every horodisk centered at $x.$ The set of all horocyclic points is denoted by $L_{h}(G).$ A point in  $L(G)$ is parabolic if there exists a parabolic transformation $g\neq Id$ in $G$ such that $g(x)=x.$ The set of all parabolic points is denoted by $L_{p}(G).$

For $g$ in $G$, we denote by $\mathcal{H}_{z}(g)$ the closed hyperbolic half-plane containing $z$, bounded by the perpendicular bisector of the segment $[z, g(z)]_{h}$. The Dirichlet fundamental domain $\mathcal{F}_{z}(G)$ of $G$ centered at $z$ is the intersection of all the sets $\mathcal{H}_{z}(g)$ with $g$ in $G-\{Id\}$. For simplicity, in this paper we use the notation $\mathcal{F}$ for the Dirichlet
fundamental domain $\mathcal{F}_{z}(G)$ of $G$ centered at $z=i.$  It is easy to see that a Dirichlet domain is a star-shaped convex subset  of $\mathbb{H}$ under the topology induced by the hyperbolic metric. For the Dirichlet fundamental domain $\mathcal{F}$, let $\mathcal{F}^{\circ}$ denote its interior and $\bar{\mathcal{F}}$ its closure. The boundary at infinity of $\mathcal{F} $ is the set denoted by $\mathcal{F}(\infty)$ defined by $\mathcal{F}(\infty)=\bar{\mathcal{F}}\cap \mathbb{R}.$ Let $ g$ be a nontrivial element of $G$. When the intersection of $\mathcal{F} $ and $g(\mathcal{F}) $ is non-empty, it is contained in  the perpendicular bisector of the segment $[z, g(z)]_{h}$. This intersection, is a point, a non-trivial geodesic segment, a geodesic ray or a geodesic. In the latter three cases, we say this intersection is an edge. The vertices are the endpoints of the edges. An infinity vertex is a vertex  contained in $\bar{\mathbb{R}}.$ When $G$ is of first kind, the  infinity vertex set $\mathcal{F}(\infty) $ is at most countable (it is possible empty). A Fuchsian group $G$ is called a lattice if, for any $z\in \mathbb{H}$, the area of each Dirichlet fundamental domain $\mathcal{F}_{z}(G)$ is finite. The readers are suggested to see \cite{Be,Da, KMS}  for more details about Fuchsian groups.
\subsection{ Carleson measures.} A positive measure $\lambda$ defined in a simply connected domain $\Omega$ of the complex plane $\mathbb{C}$ is called a Carleson measure if there exists some constant C which is independent  of $r$ such that, for all $0<r<diameter(\partial\Omega)$ and $z\in \partial\Omega$,
 $$\lambda(\Omega\cap D(z,r))\leq Cr. $$
 The infimum of all such $C$ is called the Carleson norm of $\lambda,$ denoted by
 $\parallel\lambda\parallel_{*}.$ In this paper, we mainly focus our attention on the cases $\Omega=\mathbb{H} $ or $\Omega=\mathbb{D}.$ We denote by $CM(\mathbb{H})$ ($CM(\mathbb{H})$) the set of all Carleson measures on the upper half plane $\mathbb{H}$(the unit disk $\mathbb{D}$). For more details about Carleson measure, see \cite{Ga,WZ}.

 We say that an essentially bounded measurable function $\mu(z)$ belongs to $ \mathcal{CM}(\mathbb{H})$ (or $ \mathcal{CM}(\mathbb{D})$) if it induces a Carleson measure $\lambda_{\mu}\in CM(\mathbb{H})$ (or $ \mathcal{CM}(\mathbb{D})$) by $$\lambda_{\mu}(z)=\displaystyle\frac{|\mu|^{2}}{Im(z)}dxdy\in CM(\mathbb{H})(~~\text{or}~~\displaystyle\frac{|\mu|^{2}}{1-|z|^{2}}dxdy\in CM(\mathbb{D}).$$ The importance of the class $ \mathcal{CM}(\mathbb{H})(~~\text{or}~~\mathcal{CM}(\mathbb{D}))$ lies in the fact that it plays a crucial role in the theory of BMO-{\teich} space, see \cite{AZ},\cite{Cu}, \cite{CZ},\cite{FKP} and \cite{SW}.

\subsection{Conjugations of Fuchsian groups.}
We now consider conjugations of a Fuchsian group $G$ by a quasiconformal self homeomorphisms of the upper half plane. We say that a quasiconformal  homeomorphism $f\in QC(\mathbb{H}) $  is $G$-compatible if the conjugation $f\circ g\circ f^{-1}$ is a {\mobi} transformation for every $g\in G$, denoted all such homeomorphisms by $QC(G)$. Similarly,  a quasisymmetric homeomorphism $h\in Hom^{+}(\mathbb{R})$  is called $G$-compatible if the conjugation $h\circ g\circ h^{-1}$ is a {\mobi} transformation for every $g\in G,$ denoted all such homeomorphisms by $QS(G).$ The Beltrmai coefficient (or complex dilatation) of a $G$-compatible quasiconformal mapping satisfies $\mu(z)=\mu(g(z))\overline{g'(z)}/g'(z)$ for every $g\in G$, then we say that $\mu$ is a  $G$-compatible Beltrami coefficient. We denote by $M(G)$ the set of all $G$-compatible Beltrami coefficients.

 For a $G$-compatible Beltrami coefficient $\mu$, when $\mu\in \mathcal{CM}(\mathbb{H})$ and the Carleson norm of the induced measure $\lambda_{\mu}=|\mu|^{2}(z)/Im(z)dxdy$ is small, then $f_{\mu}(\mathbb{R})$ is a locally rectifiable curve, where $f_{\mu}$ is the quasiconformal mapping of the complex plane $\mathbb{C}$ with $0$, $1$ and $\infty$ fixed, whose Beltrami coefficient equals to $\mu$ almost everywhere in $\mathbb{H}$ and equals to zero else where on the complex plane $\mathbb{C}$, see \cite{Se}. Perhaps the best way to see this is to transfer the problem from the upper plane $\mathbb{H}$ to the unit disk $\mathbb{D}$. Suppose $\mu\in \mathcal{CM}(\mathbb{D})$ and the Carleson norm of the induced measure $\lambda_{\mu}=|\mu|^{2}(z)/1-\mid z\mid^{2}dxdy$ is small, then $f_{\mu}(S^{1})$ is a chord-arc curve, where a Jordan curve $\gamma$ is called a chord-arc curve if it is  rectifiable and $L(z_{1},z_{2})\leq K|z_{1}-z_{2}|$ for any $z_{1}, z_{2}\in \gamma$, $S^{1}$ denotes the unit circle, $L(z_{1},z_{2})$ denotes the Euclidean length
of the shorter arc of $\gamma$ between $z_{1}$ and $z_{2}$ . This is essential for the proof of the convergence-type first-kind Fuchsian groups failing to have Bowen's property, see \cite{AZ1}. It is also the method to prove that some convergence-type Fuchsian groups  fail to have Ruelle's property, see \cite{HZ, HW}.
\subsection{Denjoy domain.}
Recall that a Denjoy domain is a connected open subset $\Omega$ of the extended complex plane $\bar{C}$ such that the complement $E=\bar{C}\setminus\Omega$ is a subset of the real axis $\mathbb{R}.$ In addition, $\Omega$ is called a Carleson-Denjoy domain if there exists a positive constant $C$ such that
 $$|E\bigcap(x-t,x+t)|\geq Ct $$
for all $x\in E$ and $0,t<\text{diam}(E)$, where $|\cdot|$ denotes the Lebesque measure on $\mathbb{R}$.

 We now describe a fundamental domain for the covering group of the Denjoy domains. We consider the unit disk $\mathbb{D}$ as the universal covering surface. Let $F$ be a closed subset of positive measure of the unit circle $S^{1}$. The complement of $F$ in the unit circle $S^{1}$ will be a collection $\{J_{n}\}$ of disjoint open intervals. Corresponding to each $J_{n}$, we construct the circle orthogonal to the unit circle which passes through the endpoints of $J_{n}$. Denote by $\alpha_{n}$ the arc of that circle that intersects $S^{1}$. Then $\{\alpha_{n}\}$ together with $F$ bounded a simply connected region $\Omega$ with rectifiable boundary. Let $\varphi$ be a conformal mapping of $\Omega$
 onto the upper half-plane $\mathbb{H}$. Extend $\varphi$ to the boundary of $\Omega$ and $\varphi$ preserves the sets of positive measure and sets of zero measure, since the boundary of $\Omega$ is rectifiable. We can get that the image of $F$ under $\varphi$ is just a closed subset $E$ of positive measure of $\mathbb{R},$ and each $\alpha_{n}$ will map onto a complementary open interval $I_{n}$ in $\mathbb{R}.$ A reflection across $\alpha_{n}$ yields a mapping  of the reflected region onto the lower half plane that can be extended analytically across $\alpha_{n}.$ Reflection across any $\alpha_{n}$ extends $\varphi$ to an analytic mapping of a larger region $\Omega'$ on to the complement of the closed set $E$. The region $\Omega'$ is a fundamental domain for the Denjoy domain $\mathbb{C}\setminus E,$ for more details see \cite{RR}. It is easy to see that when  $\mathbb{C}\setminus E $ is  totally disconnected, the covering group $G$ is of infinitely generated Fuchsian group of the first kind and of  convergence type.
\subsection{Main results.}The set $SQS(\mathbb{R})$ of all strongly quasisymmetric homeomorphisms of $QS(\mathbb{R})$ were much investigated because of their great important in the application to harmonic analysis, see \cite{AZ, CZ, FKP, HS, Sh, SW, WM}. The sub-class $SQS(G)$ of all strongly quasisymmetric homeomorphisms which are compatible with Fuchsian groups is also very important since it plays a very critical role to prove that some convergence-type Fuchsian groups  fail to have Ruelle's property, see \cite{AZ, HW, HZ}. There are many equivalent conditions for a quasisymmetric homeomorphism to be strongly quasisymmetric homeomorphism. In
 particular we mention here is that a quasisymmetric homeomorphism $h\in QS(\mathbb{R})$ is in $SQS(\mathbb{R})$ if and only if it has a quasiconformal extension  $f$ to the upper half plane $\mathbb{H}$ onto itself such that the induced measure $\lambda_{\mu}=|\mu|^{2}/Im(z)dxdy$ by the Beltrami coefficient $\mu$ of $f$ is a Carleson measure on $\mathbb{H}$, i.e. $\lambda_{\mu}\in CM(\mathbb{H})$ or $\mu \in \mathcal{CM(}\mathbb{H})$, see\cite{FKP }or \cite{AZ}. For a Fuchsian group $G$,  Cui and Zinsmeister showed that $h\in SQS(G)$ if and only if the Beltrami coefficient of the Douady-Earle extension of $h$ to the upper half plane is in $\mathcal{CM}(\mathbb{H})$.

It is well known that for a Fuchsian group $G$, the upper plane $\mathbb{H}$ can be tessellated  by any one Dirichlet fundamental domain $\mathcal{F}_{z}(G)$ of $G$. Is it possible  to check  that the measure $\lambda_{\mu}=|\mu|^{2}/Im(z)dxdy$ induced by $\mu\in M(G)$ is in $CM(\mathbb{D})$  directly from its value on the Dirichlet fundamental domain $\mathcal{F}_{z}(G)$ of $G$?  In this paper, we may in some sense give some positive answer to this question.

Recall that a domain $\mathcal{F}_{z}^{*}$ is called a generalized Dirichlet fundamental domain of $G$ if for each $x\in \mathcal{F}_{z}(\infty)\cap L_{p}(G)$, there exists a horodisk $B_{x}$ centered at $x$ such that $\mathcal{F}^{*} $ is the union of the Dirichlet fundamental domain  $\mathcal{F}_{z}(G)$ and the domain $\bigcup\limits_{{\rm{x\in\mathcal{F}_{z}(\infty)}}} {B_{x}}$, where $L_{p}(G)$ denotes the set of all parabolic limit points of $G$ .
In this paper, we first show
\begin{theorem}\label{main1}
Let $G$ be a convergence Fuchsian group of  first kind  and $\mathcal{F}$ the Dirichlet  fundamental domain of $G$ centered at $i$. A quasisymmetric homeomorphism $h$ is in $SQS(G)$ if and only if it has a quasiconformal extension  $f$ to the upper half plane $\mathbb{H}$ onto itself such that the Beltrami coefficient $\mu$ of $f$ is in $M(G)$ and the induced measure $\lambda_{\mu}=\frac{|\mu(z)|^{2}}{y}dxdy$ by $\mu$  is a Carleson measure on a generalized Dirichlet fundamental domain $\mathcal{F}^{*}$ of $G.$

\end{theorem}
{\bf Remark:} Notice that Theorem \ref{main1} fails for the case of the divergence Fuchsian groups , since  Bishop\cite{Bi1} showed that the divergence Fuchsian groups hold a rigidity property, now called Bowen's property, i.e. the $G$-compatible quasisymmetric homeomorphisms are always singular. Hence the measure $|\mu|^{2}/ydxdy$ can not be a Carleson measure.

This theorem means that the Carleson property of a measures which are compatible with the convergence Fuchsian groups of first kind can be checked from its value  in  a generalized Dirichlet fundamental domain $\mathcal{F}^{*}$.

When a Fuchsian group $G$ is of first kind, $\mathcal{F}(\infty)$ contains at most countable many points.  Let $G$ be a Fuchsian group of second kind. Then $\mathcal{F}(\infty)$ contains at least one non-trivial interval $J$ whose endpoints are two distinct infinite vertices of $\mathcal{F}.$  First of all we recall a result in \cite{H1} which considered the case of  the finite generated Fuchsian groups of second kind acting on the unit disk $\mathbb{D}=\{z\in \mathbb{C}:|z|<1\}$ as follows.

\begin{theorem}\label{gs1}\cite{H1}
Let $G$ be a finite generated Fuchsian group of the second kind and $\mathcal{F}$ the Dirichlet  fundamental domain of $G$ centered at $0$. Suppose $\mu$ is in $ M(G).$  The induced measure $(|\mu|^{2}/{1-|z|^{2}})dxdy$ is in $ CM(\mathcal{F})$ if and only if
$(|\mu|^{2}/{1-|z|^{2}})dxdy$ is in $ CM(\mathbb{D}).$
\end{theorem}
A modest motivation of this paper was to make a correction of a mistake about Theorem \ref{gs1}  which holds for all finite generated Fuchsian group of second kind. In this theorem, we cannot drop the assumption that the finite generated Fuchsian group $G$ does not contain the  parabolic elements.  For example when $G$ is a cyclic parabolic group, e.g., the covering group of a punctured disk. In this case the fundamental domain in the upper half-plane is a vertical strip and setting $|\mu|$ to be constant in the part of the strip above height 1  to provide a counterexample, since if $|\mu|$ is constant on a horodisk, then it does not give a Carleson measure.
In this paper, we revise Theorem \ref{gs1} as
\begin{theorem}\label{main2}
Let $G$ be a finite generated Fuchsian group of the second kind. Let $\mathcal{F}$ be the Dirichlet  fundamental domain of $G$ centered at $0$.  Suppose $\mu$ is in $ M(G).$ Then $\mu$ is in $\in\mathcal{ CM}(\mathbb{D})$ if and only if there exists a generalized Dirichlet fundamental domain $\mathcal{F}^{*}$ of $G$ such that  the induced measure $(|\mu|^{2}/{1-|z|^{2}})dxdy$ by $\mu$ is a Carleson measure in $ \mathcal{F}^{*}.$
\end{theorem}

Recall that a Denjoy domain is a connected open subset $\mathcal{D}$ of the extended complex plane $\bar{C}$ such that the complement $E=\bar{C}\setminus\mathcal{D}$ is a subset $E$ of the real axis $\mathbb{R}.$ In this paper we consider the case that $E$ is the union of countable many non-trivial intervals $\{I_{n}\}.$ For such domains, we can generalize Theorem \ref{gs1} to infinitely generated Fuchsian groups as

 \begin{theorem}\label{main3}
Let $\mathcal{D}$ be any Carleson-Denjoy domain with at most countable many non-trivial closed interval boundary(may contain infinity point $\infty$). Let $G$ be the covering group of $\mathcal{D}$ with the unit disk as the universal covering surface. Let $\mathcal{F}$ be a  fundamental domain of $G$ and $\mu\in M(G)$. Then $\mu\in \mathcal{CM}(\mathbb{D})$  if and only if the measure
$(|\mu|^{2}/1-|z|^{2})dxdy$ is  a Carleson measure in the fundamental domain $\mathcal{F}.$
\end{theorem}

We will consider another extremity case of Denjoy domains that the similar result does not hold as Theorem \ref{main3} . We shall use the rich knowledge about the Ruelle's property about Fuchsian groups to prove the following  result.

\begin{theorem}\label{main4}
  Suppose  $(s_n)$ is a sequence of real numbers increasing to infinity and $G$ is the covering group of the surface $S=\mathbb{C}\backslash\{s_n,\,n\geq 0\}.$  There exists a sequence $(s_n)$ such that  any $G$-compatible quasi-symmetric homeomorphisms are always singular.
\end{theorem}

\subsection{Notation.} The expressions  $A \lesssim B$, $A\gtrsim B$ and $A\asymp B$ mean that $A<C B$, $A>C B$ and $B/C<A<CB$  for an unspecified constant $C.$

\section{Some lemmas}

 Let us recall that a Fuchsian group $G$ has Ruelle's property if, for any  family of Beltrami coefficients $(\mu_{t})\in M(G)$ which is analytic in $t\in\Delta$, the map $t\mapsto HD(L(G_{\mu_{t}}))$ is  real-analytic in $\Delta$, where $HD(L(G_{\mu_{t}}))$ denotes the Hausdorff dimension of the the limit set of the quasi-Fuchsian groups $G_{\mu_{t}}$ .
In 1982, Ruelle \cite{Ru} showed that all cocompact groups have this property. In 1997, Anderson and  Rocha \cite{AR} extended this result to finitely-generated Fuchsian groups without parabolic elements.  In \cite{AZ, AZ1}, Astala and Zinsmeister showed that for Fuchsian groups corresponding to Denjoy-Carleson domains or infinite $d$-dimensional "jungle gym" with $d\geq 3$ , Ruelle's property fails.  It is easy to see that the Fuchsian groups
studied in \cite{AZ1}and \cite{Bi1} are of first kind.
Very recently, the first author and Zinsmeister showed that

\begin{lemma}{\label{le1}}\cite{HZ}
All convergence type Fuchsian groups of the first kind fail to have Ruelle's property.
  \end{lemma}

  In \cite{HZ}, we constructed an infinitely generated Fuchsian group which has Ruelle's property. As the author's known, this is the first concrete example of infinitely generated Fuchsian group which has Ruelle's property.

\begin{lemma}\label{le2}\cite{HZ}
 There exists a sequence $(s_n)$ of real numbers increasing to infinity such that the Fuchsian group uniformizing $S=\mathbb{C}\backslash\{s_n,\,n\geq 0\}$ has Ruelle's property.
\end{lemma}

In this paper, we will  need the following lemma which
 essentially belongs to Astala and Zinsmeister, see \cite{AZ}, \cite{AZ1}, or (\cite{H1}, Lemma 2.2).
\begin{lemma}\label{le3}
For a convergence-type Fuchsian group $G$ acting on the unit disk $\mathbb{D}$ and $\mu$ in $M(G)$, if the support set of $\mu\chi_{\mathcal{F}}$ is contained in a compact subset of the unit disk $\mathbb{D}$, then $\mu$ is in $\mathcal{CM}(\mathbb{D}).$
\end{lemma}
A result from \cite{Z} says that

\begin{lemma}\label{le4}(\cite{Z})
Let $\Omega$ be a chord-arc domain. Then the following are equivalent:

(a) $d\nu$ is a Carleson measure for $\Omega.$

(b) For $0<p<\infty$ and $f\in H^{p}(\Omega),$
$$\iint_{\Omega}|f|^{p}dv\leq C\int_{\partial \Omega}|f|^{p}ds,$$
where $H^{p}(\Omega)=\{f: f \text{is analytic on }\Omega \text{ and } \int_{\partial \Omega}|f|^{p}ds<\infty \}$ and the constant $C$ depends only on the the Carleson norm of $d\nu.$
\end{lemma}

\section{Infinity boundary of Dirichlet domains}
Let $G$ be a Fuchsian group acting on the upper half plane. The goal of this section is to give characterizations of the infinity boundary of a Dirichlet fundamental domain $\mathcal{F}$ of $G$ in terms of limit points. In this paper, we mainly consider the case of infinitely generated Fuchsian groups, the finite generated case can see \cite{Fr}.

Suppose $G$ is Fuchsian group of first kind. If the parabolic limit set $L_{p}(G)$ of $G$ is not the empty set. Then $L_{p}(G)\cap \mathcal{F}(\infty) $ is not the empty set. Moreover each point in  $L_{p}(G)\cap \mathcal{F}(\infty) $  is isolated and is an infinity vertex of  $\mathcal{F}.$  For $x\in L_{p}(G)$, denote the elements in $G$ fixing $x$ by $G_{x}$. It is known that $G_{x}$ is a cyclic subgroup of $G.$ Suppose the group $G_{x}$ is generated by a parabolic isometry and hence preserves each horodisk $B_{x}$ centered at $x.$ The quotient surface $B_{x}/G_{x} $ is called a cusp associated with $B_{x}.$

It is known that when $G$ is a non-elementary finitely generated Fuchsian group of first kind. The infinity boundary $\mathcal{F}(\infty)$ is finite, and is equal to the set $L_{p}(G)\cap \mathcal{F}(\infty) $. For infinitely generated Fuchsian groups,  we need the following
\begin{lemma} (\cite{Fr},Page 39)\label{fr3}
If $G$ is a infinitely generated Fuchsian group, for any Dirichlet fundamental domain $\mathcal{F}_{z}(\infty)$, there exists at least one point in $\mathcal{F}_{z}(\infty)$ which is not in $L_{p}(G).$
\end{lemma}
To better understand the properties of the infinity boundary $\mathcal{F}(\infty)$ of the Dirichlet fundamental domain $\mathcal{F}_{z}$, we need a binary classification of the limit points.  A point $x\in \bar{\mathbb{R}}$ is a conical point or radial point of $G$ if there is a sequence of elements $g_{i}\in G$ such that for any $z\in \mathbb{H}$, there exists a constant $C$ and a hyperbolic line $\gamma$ with endpoint $x$ such that the hyperbolic distance between $g_{i}(z)$ and $\gamma$ are bounded by
$C.$ Denote by $L_{c}(G)$ the set of all the conical limit points and $L_{e}(G)$ the set of all the escaping limit points. Let $S=\mathbb{H}/G$ be the corresponding surface of $G$. The points in $L_{c}(G)$ are just corresponding to the geodesics in $S$ which return to some compact set infinitely often and the points in $L_{e}(G)$ are corresponding to the geodesics in $S$ which eventually leave every compact subset of $S$. When $G$ is an infinitely-generated first-kind Fuchsian group, the set $\mathcal{F}(\infty)$ is contained in $L_{e}(G).$ We first make clear  the relationships between the conical limit points set $L_{c}(G)$, the escaping limit points set $L_{e}(G)$, the parabolic limit points set $L_{p}(G)$  and the horocyclic limit points set $L_{h}(G)$, respectively. It is known that $L_{p}(G)$ is contained in $L_{e}(G)$, $L_{c}(G)$ is contained in $L_{h}(G)$ and a parabolic limit point is not horocyclic. In this paper, we need the following lemmas.

\begin{lemma}(\cite{Fr}, Page 30)\label{fr1}
Let $G$ be a non-elementaty Fuchsian group and $x$ be a point in  $L_{p}(G)$. There exists a horodisk $B_{x}$ centered at $x$ such that $g(B_{x})\cap B_{x}=\emptyset,$ for any $g$ in $G$ which does not fix $x.$
\end{lemma}
\begin{lemma}(\cite{Fr}, Page 36)\label{fr2}
Let $G$ be a non-elementary Fuchsian group. Then the set  $L_{h}(G)\cap \mathcal{F}(\infty) $ is the empty set.
\end{lemma}

Now we give the proofs of Theorem \ref{main1}.

\section{Proof of Theorem \ref{main1}}
A quasisymmetric homeomorphism $h\in QS(G)$ is in $SQS(G)$ if and only if it has a quasiconformal extension  $f$ to the upper half plane $\mathbb{H}$ onto itself such that the induced measure $\lambda_{\mu}=|\mu|^{2}/Im(z)dxdy$ by the Beltrami coefficient $\mu$ of $f$ is a Carleson measure on $\mathbb{H}$, i.e. $\lambda_{\mu}\in CM(\mathbb{H})$ or $\mu \in \mathcal{CM}(\mathbb{H})$, see\cite{FKP }or \cite{AZ}. Hence if $h\in SQS(G)$, the proof of the existence of a generalized Dirichlet fundamental domain $\mathcal{F}^{*}$ such that $\lambda_{\mu}\in CM(\mathcal{F}^{*})$ is trivial.  We only need to prove that if there exists a generalized Dirichlet fundamental domain $\mathcal{F}^{*}$ such that  a quasisymmetric homeomorphism $h\in QS(G)$ has a quasiconformal extension $f$ to the upper half plane $\mathbb{H}$ onto itself with Beltrami coefficient $\mu$ satisfying $(|\mu|^{2}/y)dxdy \in CM(\mathcal{F}^{*}),$ then $(|\mu|^{2}/y)dxdy \in CM(\mathbb{H})$.

Since Carleson measure property is preserving by conformal map between simple connected domains, in order to avoid analysing the infinity point $\infty$ alone, here we use the unit disk $\mathbb{D}$ as the universal covering surface.

Suppose $G$ is a convergence type Fuchsian group of first kind, it must be infinitely generated.  We first show the following claim.

 \noindent{\bf Claim 1:} There exists a generalized Dirichlet  fundamental domain $\mathcal{F}^{*}$ such that the horodisks in the set

 $$\{g(B_{x}): g\in G, x\in L_{p}(G)\cap\mathcal{F}(\infty)\}$$
 are pairwise either disjoint or identical.

 Note that each point in $L_{p}(G)\cap \mathcal{F}(\infty)$ is isolated in $\mathcal{F}(\infty)$, we can choose  the horodisks centered at points in $L_{p}(G)\cap \mathcal{F}(\infty)$ sufficiently small such that they are pairwise disjoint.
 Moreover by Lemma \ref{fr1}, we associate to each point $x$ in $L_{p}(G)\cap \mathcal{F}(\infty)$ a horodisk $B_{x}$ centered at $x$ satisfying the condition of
 $$B_{x}\cap g(B_{x})=\emptyset$$
  for all $g\in G$ not fixing $x.$

 If the set $L_{p}(G)\cap \mathcal{F}(\infty)$ contains at most one point, by Lemma \ref{fr1}, the claim holds trivially.
Otherwise, the set $L_{p}(G)\cap \mathcal{F}(\infty)$ contains at least two points (maybe infinity many points).

Suppose the set $L_{p}(G)\cap \mathcal{F}(\infty)$ contains finite many points. Let $x$ and $y$ be in $L_{p}(G)\cap \mathcal{F}(\infty)$.
 If $y$ is in the orbit $Gx$. Then there exists $g\in G$ such that $g(B_{x})$ is contained in $B_{y}$ or $B_{y}$ is contained in $g(B_{x})$.
  Hence we can let $B_{y}=g(B_{x})$ in the former case or $B_{x}=g^{-1}(B_{y})$ in the later case.

  We suppose $y$ is not in the orbit $Gx.$ Consider the set $A$ of $g\in G$ such that $g(B_{x})\cap B_{y}\neq\emptyset.$ Let $Q$ denote the set of two sided cosets $G_{y}\setminus A/G_{x}$ (i.e., equivalence classes via two relations). If the set
  $Q$ is a finite set, hence we can replace $B_{x}$ with a smaller horodisk than $B_{x}$ such that $g(B_{x})\cap B_{y}$ is empty set for all $g\in G.$ Otherwise, there are infinity many $g\in G$, such that $g(B_{x})\cap B_{y}\neq\emptyset.$ Consider a sequence of distinct elements of $Q$ written as $(q_{n}=G_{y}a_{n}G_{x})_{n\geq 1}$ where $a_{n}$ is an element in $A$. Fix a compact fundamental domain $K$ for the action of $G_{y}$ on the horodisk $B_{y}$.
 We have for each positive integer $n$, there exists $g_{n}\in G_{y}$ such that $g_{n}\circ a_{n}(B_{x})\cap K$  is not empty. On the other side, for fixed $x\in L_{p}(G)\cap\mathcal{ F}(\infty)$, since the choosing $B_{x}$ is by Lemma \ref{fr1} which are pairwise disjoint, this is contradict with $g_{n}\circ a_{n}(B_{x})\cap K\neq\emptyset.$ Hence for each pair points $x, ~y\in L_{p}(G)\cap \mathcal{F}(\infty),$ there exists the horodisks $B_{x}$ and $B_{y}$ such that the images of them are pairwise either disjoint or identical under the transmission of $G.$

 Suppose the set $L_{p}(G)\cap \mathcal{F}(\infty)$ contains infinite many points.
 By the proof of finite points case and applying Mathematical induction method, we can easily show that the claim holds.

By Lemma \ref{fr3}, we know that the set $\mathcal{ F}(\infty)\setminus L_{p}(G)$ is not empty set. By  Lemma \ref{fr2}, we know that the points $x\in \mathcal{ F}(\infty)\setminus L_{p}(G)$ are not horocyclic. Furthermore for points in  $x\in \mathcal{ F}(\infty)\setminus L_{p}(G)$, we have

\noindent{\bf Claim 2:} There exists a horodisk $B_{x}$ centered at $x\in \mathcal{ F}(\infty)\setminus L_{p}(G)$ such that $B_{x}\subset\mathcal{ F}.$

Fix $x\in \mathcal{ F}(\infty)\setminus L_{p}(G)$, to show this claim, we only need to show that there exists a horodisk $B_{x}$ such that for any non-identity element $g\in G$, $g(B_{x})$is disjoint with $ B_{x}$.
Otherwise, suppose that for any horodisk $B_{x}$, there exists $g\in G$ and $g\neq Id$ such that $g(B_{x})\cap B_{x}=\emptyset.$ We can choose a sequence $(B_{n})_{n\geq1}$ of horodisks centered at $x$  with radii $\displaystyle\frac{1}{n}$  and a sequence $(g_{n})_{n\geq 1}\in G$  with $g_{n}(B_{n})\cap B_{n}\neq\emptyset.$ Note that $G$ is discrete, hence for sufficiently large $n$, we have $g_{n}(x)=x$, this is impossible since $x$ is not in $L_{p}(G).$ Hence Claim 2 holds.

For each $x\in \mathcal{ F}(\infty)\setminus L_{p}(G) $, we can replace each horodisk $B_{x}$ with a larger domain $V_{x}$ satisfying $|V_{x}|\asymp |B_{x}|$ and $V_{x}\subset \mathcal{F}$ such that
$\mathcal{F}$ minuses the union of $$ \bigcup\limits_{{\rm{x \in {L_p}(G)\cap \mathcal{F}(\infty)} }} {{B_x}} \text{~~and~~}  \bigcup\limits_{{\rm{x \in { \mathcal{F}(\infty)\setminus L_p}(G)} }} {{V_x}}$$ being a compact subset of $\mathbb{D},$ denoted it by $K.$
Furthermore we denote $$W= \bigcup\limits_{{\rm{x \in {L_p}(G)\cap \mathcal{F}(\infty)} }} {{B_x}} \text{~~and~~} V= \bigcup\limits_{{\rm{x \in { \mathcal{F}(\infty)\cap L_p}(G)} }} {{V_x}},$$ respectively.
Then the generalized Dirichlet fundamental domain $\mathcal{F}^{*}$  divided into three pairwise disjoint parts(except a set of zero measure) as $K$, $W$ and $V$, and the unit disk $\mathbb{D}$ is divided into three pairwise disjoint parts (except a set of zero measure) as $ \bigcup\limits_{{\rm{g \in G} }} {{g(W)}}$, $  \bigcup\limits_{{\rm{g \in G} }} {{g(V)}}$ and $ \bigcup\limits_{{\rm{g \in G} }} {{g(K)}}$, denoted by $\Omega_{1}$, $\Omega_{2}$ and $\Omega_{3}$, respectively.

Now its time to show that if $$\displaystyle \frac{|\mu|^{2}\chi_{\mathcal{F}^{*}}}{1-|z|^{2}}\in CM(\mathcal{F}^{*})\text{~~ then~~ }\displaystyle \frac{|\mu|^{2}}{1-|z|^{2}}\in CM(\mathbb{D}).$$

For a Fuchsian group $G$ of convergence type and  any $\mu\in M(G)$, note that $K$ is a compact subset of the Dirichlet fundamental domain $\mathcal{F}$, thanks to Lemma \ref{le3} we have the measure $$\displaystyle \frac{|\mu|^{2}\chi_{\Omega_{3}}}{1-|z|^{2}}\in CM(\mathbb{D}).$$

In the following we only need to show that both the measures $$\displaystyle \frac{|\mu|^{2}\chi_{\Omega_{1}}}{1-|z|^{2}}  \text{~~and~~} \displaystyle \frac{|\mu|^{2}\chi_{\Omega_{2}}}{1-|z|^{2}}$$  are in $CM(\mathbb{D})$.

Now we need an estimate about how close the disjoint horodisks can pack around the base of a given one horodisk $B^{*}$.  Suppose $B_{n}$ is a horodisk which is based at points (the tangency point of $B_{n}$ with the unit circle $S^{1}$ ) whose Euclidean distance from the center of $B^{*}$  lies between $\frac{1}{n}$ and $\frac{1}{n+1}$.

 \noindent{\bf Claim 3:} There exist a constant $N$ depending only on the radius of $B^{*}$ such that  for $n>N$ , $B_{n}$ has radius $\lesssim \frac{1}{n^{2}}.$

 This claim is from \cite{Su } where is for the case of Kleinian groups.

The total available area under $B^{*}$ is that of a region whose width is about $(\frac{1}{n}-\frac{1}{n+1})+\frac{2}{n^{2}}$ and whose height is about $(\frac{1}{n}+\frac{1}{n^{2}})^{2}.$ Since each of these has order $\frac{1}{n^{2}}$, the area has order $\frac{1}{n^{4}}.$

Let $x$ be a point in $\mathcal{ F}(\infty)\cap L_{p}(G)$. Let $\mathcal{G}_{x}$ denote the coset space $G/G_{x}$. It is easy to see that the images of the horodisk $B_{x}$ under the translation of $G$ only depends on the coset of $g\in \mathcal{G}_{x}$. The measure of the image $g(B_{x})$ of the horodisk $B_{x}$  under the measure $\frac{{{{\left| {{\mu(w)}} \right|}^2}}}{{1 - \left| w \right|^{2}}}dudv$ is
\begin{eqnarray}
&&Measure(g(B_{x}))=\iint_{g(B_{x})}\frac{{{{\left| {{\mu(w)}} \right|}^2}}}{{1 - \left| w \right|^{2}}}dudv\\
&=& \iint_{B_{x}}\frac{{{{\left| {{\mu(g(z))}} \right|}^2}}}{{1 - \left| g(z) \right|^{2}}}\mid g'(z)\mid^{2} dxdy\\
&=& \iint_{B_{x}}\frac{{{{\left| {{\mu(g(z))}}\frac{\overline{g'(z)}}{g'(z)} \right|}^2}}}{{1 - \left| g(z) \right|^{2}}}\mid g'(z)\mid^{2} dxdy\\
&=&\iint_{B_{x}}\frac{\mid\mu(z)\mid^{2}}{1-\left|z\right|^{2}}\left|g'(z)\right|dxdy\\
&\lesssim&\int_{\partial B_{x}}\left|g'(z)\right|ds
=\int_{\partial g(B_{x})}ds= 2\pi r_{g(B_{x})},
\end{eqnarray}
where the equality (5.4) holds by $\mu\in M(G)$ and the inequality $\lesssim$ in (5.5)   is by Lemma \ref{le4} for the domain $B_{x}$ to $g'$ for the case $p=1.$
Thanks to Claim 1, we know that for $x, y\in \mathcal{ F}(\infty)\cap L_{p}(G)$ and
$g_{1}\in \mathcal{G}_{x}$ and $g_{2}\in \mathcal{G}_{y}$, $g_{1}(B_{x})\cap g_{2}(B_{y})=\emptyset.$

To show the measure $\displaystyle \frac{|\mu|^{2}\chi_{\Omega_{1}}}{1-|z|^{2}}$ is in $CM(\mathbb{D})$, we use the  interpolating sequences: Recall that a sequence $\{z_{j}\}$ is called an interpolating sequence of $\mathbb{D}$ if $$(i) ~~\exists \delta>0, ~\rho(z_{j}, z_{k})\geq\delta~~if ~~j\neq k; $$
 $$(ii)~~\sum (1-|z_{i}|^{2})\delta_{z_{i}}\in CM(\Delta),$$
 where $\delta_{z}$ stands for the Dirac mass at $z$.

The sequence $\{g(0)\}_{g\in G}$
is an interpolating sequence of the unit disk $\Delta$ when $G$ is convergence type, see \cite{AZ} or for the detail of the proof see \cite{H1}.
 Note that $x\in \mathcal{ F}(\infty)\cap L_{p}(G)$ and $g\in G/G_{x}$, the Euclidean radii of the horodisks centered at $x$ and $g(x)$ are in the ratio $|g'(x)|$ up to a universal factor(see \cite{Su}, page 291). Hence for  $x\in \mathcal{ F}(\infty)\cap L_{p}(G)$, we can choose $\xi\in \partial B_{x}\cap \mathcal{F}$ such that
 $$1-g(\xi)\asymp r_{g(B_{x})}\asymp |g'(x)|(1-|\xi|)\asymp |g'(x)|r_{B_{x}},$$
 where $r_{B_{x}}$ denotes the radius of the horodisk $B_{x}.$
 Furthermore, the inequality in (5.5) holds only depending on the Carleson norm of the meaure $(|\mu|^{2}/1-|z|^{2})dxdy$ on the generalized Dirichlet fundamental domain $\mathcal{F}^{*}$ and $$ \Omega_{1}=\bigcup\limits_{{\rm{x \in {L_p}(G)\cap \mathcal{F}(\infty)} }} {\bigcup\limits_{{\rm{g \in  \mathcal{G}_{x}}}} {g(B_x)}},$$
the measure $\displaystyle \frac{|\mu|^{2}\chi_{\Omega_{1}}}{1-|z|^{2}}\in CM(\mathbb{D})$ is proved.

It is known that when the generators of $G$ contain parabolic elements, the boundary of the Dirichlet fundamental domain is not a chord-arc curve since it contains cusp points. Hence the using of Lemma \ref{le4} may be fail.
However the boundary of $\mathcal{F}\setminus W$ does not contain cusps any more. By Claim 2 and Claim 3, for each $x\in \mathcal{F}(\infty)\setminus L_{p}(G)$, we can
choose a domain $V_{x}$ such that $B_{x}\subset V_{x}\subset \mathcal{F}$ and $|V_{x}|\asymp |B_{x}|$ , and $Meas(V_{x})\asymp Meas(B_{x})$ under the measure $(|\mu|^{2}/1-|z|^{2})dxdy$ of the unit disk $\mathbb{D}$. Furthermore, we can assume that each connected component of $V$ is a chord-arc domain. It is obvious that the images of $V$ under the action of $G$ are pairwise disjoint since $V$ is a sub-domain of $\mathcal{F}$. As a copy of the proof of the measure $|\mu|^{2}\chi_{\Omega_{1}}/(1-|z|^{2})dxdy\in CM(\mathbb{D})$, we know that the measure $$|\mu|^{2}\chi_{\Omega_{2}}/(1-|z|^{2})dxdy$$ is a Carleson measure of $\mathbb{D}.$

Combine above results, the measure $$\displaystyle\frac{|\mu|^{2}}{1-|z|^{2}}dxdy=(\displaystyle
\frac{|\mu|^{2}\chi_{\Omega_{1}}}{1-|z|^{2}}+\displaystyle
\frac{|\mu|^{2}\chi_{\Omega_{2}}}{1-|z|^{2}}+\displaystyle
\frac{|\mu|^{2}\chi_{\Omega_{3}}}{1-|z|^{2}})dxdy$$ is a Carleson measure on the unit disk. Hence the measure $(|\mu|^{2}/y)dxdy$ as a pull back measure of the upper Half plane $\mathbb{H}$ by Cayley transformation $\gamma(z)=\displaystyle\frac{z-i}{z+i}$
is  a Carleson measure on $\mathbb{H}$.

\section{Proof of Theorem \ref{main2}}
Let $G$ be a finitely generated Fuchsian group of second kind. Note that for a finitely generated Fuchsian group case, by the Margulis' decomposition of  an arbitrary hyperbolic Riemann surface $\mathbb{H}/\Gamma,$  the inverse images in the upper half plane $\mathbb{H}$ of the various cusps in $\mathbb{H}/\Gamma$ consist of a  collection of disjoint horodisks. Claim 1 in the proof of Theorem \ref{main1}  naturally holds.

The second difference between Theorem \ref{main1} and Theorem \ref{main2} is that when $G$ is second kind, the boundary at infinity of  a Dirichlet fundamental domain $\mathcal{F}_{z}(G)$ contains several closed intervals whose endpoints are two distinct infinite vertices. Consider the finite sequence $(J_{i})_{1\leq i\leq k}$ of all such intervals. These intervals are pairwise disjoint. We can construct a polygon domain bounded by $J_{i}$ and the geodesics $L_{i}$ whose endpoints are on the edges of $\mathcal{F}_{z}(G)$ which have one endpoints the same as $J_{i}.$ By the same discussion as the proof of Theorem \ref{main1},this theorem can easily prove to be true, or see (\cite{H1}, Theorem 1.1) for the proof of Fuchsian group without parabolic elements.

\section{Proof of Theorem \ref{main3}}
We only need to show that the equality
 $$\iint_{B(\xi,r)\cap \mathbb{D}}\frac{{{{\left| {{\mu(w)}} \right|}^2}}}{{1 - \left| w \right|^{2}}}dudv\lesssim r,$$
 holds for all $\xi$ in the unit circle $S^{1}$ and all $0<r<2,$
  where $B(\xi, r)$ is the disk with center $\xi$ and radius $r$.

  Note that
\begin{eqnarray}
&&\iint_{B(\xi,r)\cap \mathbb{D}}\frac{{{{\left| {{\mu(w)}} \right|}^2}}}{{1 - \left| w \right|^{2}}}dudv\\
&=& \iint_{g^{-1}(B(\xi,r)\cap \mathbb{D})}\frac{{{{\sum_{g\in G}\left| {{\mu(g(z))}} \right|}^2\chi_{g(\mathcal{F})}}}}{{1 - \left| g(z) \right|^{2}}} dxdy\\
&=& \sum_{g\in G}\iint_{g^{-1}(B(\xi,r)\cap \mathbb{D})\cap \mathcal{F}}\frac{{{{\left| {{\mu(g(z))}}\right|}^2}}| g'(z)|^{2}}{{1 - \left| g(z) \right|^{2}}} dxdy\\
&=&\sum_{g\in G}\iint_{g^{-1}(B(\xi,r)\cap \mathbb{D})\cap \mathcal{F}}\frac{\mid\mu(z)\mid^{2}}{1-\left|z\right|^{2}}\left|g'(z)\right|dxdy
\end{eqnarray}

The equality $(4.3)$ holds since $\mu $ is compatible with the group $G.$
By the statement of the theorem, the infinity boundary of the fundamental domain $ \mathcal{F}(\infty)$ is the union of at most countable non-trivial closed interval, we know that the domain $\mathcal{F}(\infty)$ is a chord-arc domain. Furthermore, the boundary of the domain $g^{-1}(B(\xi,r)\cap \mathbb{D})$ is the union of two sub-circle. Hence every connected component of $g^{-1}(B(\xi,r)\cap \mathbb{D})\cap \mathcal{F} $ is a chord-arc domain.

Since the measure
$$\displaystyle\frac{|\mu|^{2}\chi_{\mathcal{F}}}{1-|z|^{2}}dxdy $$
is a Carleson measure on the domain $\mathcal{F}$, hence the measure
$$\displaystyle\frac{|\mu|^{2}\chi_{\mathcal{F}}}{1-|z|^{2}}dxdy $$ is still a Carleson measure on every component of  $g^{-1}(B(\xi,r)\cap \mathbb{D})\cap \mathcal{F} $, where $\chi_{\mathcal{F}}$ is the characteristic function of $\mathcal{F}$ .

Using Lemma \ref{le4} for the case $p=1$, we have

\begin{eqnarray}
\iint_{B(\xi,r)\cap \mathbb{D}}\frac{{{{\left| {{\mu(w)}} \right|}^2}}}{{1 - \left| w \right|^{2}}}dudv&\lesssim&\sum_{g\in G}\int_{\partial(g^{-1}(B(\xi,r))\cap \mathcal{F})}\left|g'(z)\right|ds\\
&=&\sum_{g\in G}\int_{\partial(B(\xi,r)\cap g(\mathcal{F}))}ds
\end{eqnarray}

 For Carleson-Denjoy domain $\mathcal{D}$ ,  Fernandez and Hamilton \cite{FH} showed that
\begin{eqnarray}
\sum_{\gamma\in \Gamma}\text{length}(\partial(\gamma(\mathcal{F})))<\infty,
\end{eqnarray}
 or see \cite{F}. In fact , Carleson \cite{C} for Carleson-Denjoy domain (in \cite{C} it is called homogeneous), the harmonic measure is  absolutely continuous.

  Note that $$\text{length}~~\partial (B(\xi,r)\cap g(\mathcal{F}))\lesssim ~\text{length}~(B(\xi,r)\cap\partial g(\mathcal{F})),$$
  combine with (3.6) we prove the theorem.

\section{Proof of Theorem \ref{main4}}
 We use the upper half plane $\mathbb{H}=\{z: Im(z)>0\}$ as the covering group of the surface $S.$
Let $\mathcal{D}^{*}_{1}$ be the closed disk with diameter $[0, 2]$ and $\mathcal{D}^{*}_{n}$, $n\geq 2$,  the closed disk with diameter $[2^{n-1}, 2^{n}].$
We consider the domain $$\Omega=\mathbb{H}\setminus((\cup_{n\geq1}\mathcal{D}^{*}_{n})
\cup(\cup_{n\geq 1}(-\mathcal{D}^{*}_{n}))).$$

Let $\phi$ be the conformal mapping from $\Omega$ onto $\mathbb{H}$ fixing $0$, $2$ and $\infty.$ We put $z_{0}=0$, and $z_{n}=\phi(2^{n}), \, n\geq1$ and  $z_{n}=\phi(-2^{-n}), \, n\leq-1.$ Let $\sigma_{n}$ be the reflection with respect to $\partial\mathcal{D}^{*}_{n}$ and $\tau(z)=-\bar{z}$.
By Rubel and Ryff's construction \cite{RR} of the covering group of Riemann surface
$S=\mathbb{C}\setminus \{z_{n}\},$ the Fuchsian group $\Gamma$ generated by $\{\tau\circ\sigma_{n}\}^{\infty}_{n=1}$ uniformities the surface $S$, in the sense that $S\simeq \mathbb{H}/\Gamma.$ By the construction of $\Omega$ we can see that $\Gamma$ is of infinitely generated  and of first kind. It is easy to see that $\Gamma$ contains infinitely many parabolic elements. Hence the Dirichlet domain $\mathcal{F}$ of $\Gamma$ contains countably many cusps and the set $\mathcal{F}(\infty)$, i.e., the infinity boundary  of $\Gamma$, contains countably many points, denoted  $\mathcal{F}(\infty)$ by $\{\zeta_{n}\}^{\infty}_{-\infty}.$

In the following we show that for any $\mu\in M(\Gamma)$, $\lambda_{\mu}$ is not in $ CM(\mathbb{H}).$

 By the construction of the Denjoy domain $\Omega$ and Lemma \ref{le3} we know that the group $\Gamma$ has Ruelle's property. By Lemma \ref{le2} we have that $\Gamma$ is of divergence type.  Hence the quasisymmetric homeomorphisms in $SQS(\mathbb{R})$ are singular. Note that a quasisymmetric homeomorphism $h\in QS(\mathbb{R})$ is in $SQS(\mathbb{R})$ if and only if it has a quasiconformal extension  $f$ to the upper half plane $\mathbb{H}$ onto itself such that the induced measure $\lambda_{\mu}=|\mu|^{2}/Im(z)dxdy$ by the Beltrami coefficient $\mu$ of $f$ is a Carleson measure on $\mathbb{H}$, i.e. $\lambda_{\mu}\in CM(\mathbb{H})$ or $\mu \in \mathcal{CM}(\mathbb{H})$, see\cite{FKP }or \cite{AZ}. Hence for any non-trivial(i.e. the essential norm $\parallel\cdot\parallel_{\infty}\neq0$) $\Gamma$-compatible Beltrami coefficient $\mu$ and any generalized fundamental domain $\mathcal{F}^{*}$, we always have that even though $\mu$ is in $\mathcal{CM}(\mathcal{F}^{*})$,  $\mu$ is not in $\mathcal{CM}(\mathbb{H})$.

\section{ Acknowledgements}

The first author is pleasure to thank professor Michel Zinsmeister for inviting him to the University of Orl\'{e}ans as a visiting scholar for one year and for some discussions on topics related to this paper.

\end{document}